\documentclass{elsart3-1}

\usepackage{amssymb}

\usepackage[english,francais]{babel}

\newtheorem{theorem}{Theorem}[section]

\newtheorem{e-proposition}[theorem]{Proposition}

\newtheorem{e-definition}[theorem]{Definition\rm}

\input Bach.mac
\newcommand{\Q}{\mathbb{Q}}
\newcommand{\alp}{\alpha}
\newcommand{\bet}{\beta}
\renewcommand{\ggg}{\mathfrak{g}}

\begin{document}

\centerline{}
\begin{frontmatter}


\selectlanguage{english}
\title{\huge Bach-flat Lie groups in dimension $4$}

\selectlanguage{english}
\author[label1]{Elsa Abbena},
\ead{elsa.abbena@unito.it}
\author[label1]{Sergio Garbiero}
\ead{sergio.garbiero@unito.it}
\author[label2]{\kern-4pt, Simon Salamon}
\ead{simon.salamon@kcl.ac.uk}

\vskip20pt

\address[label1]{Departimento di Matematica, Universit\`a di Torino, Via Carlo Alberto 10, 10123 Torino, Italia}
\address[label2]{Department of Mathematics, King's College London., Strand, London WC2L 2RS, London, UK}

\begin{abstract}
{\it We establish the existence of solvable Lie groups of dimension $4$ and left-invariant Riemannian
  metrics with zero Bach tensor which are neither conformally Einstein nor half conformally flat.}
\end{abstract}
\end{frontmatter}

\vskip0.5\baselineskip

\selectlanguage{english}

\section{Introduction}

Let $(M,g)$ be a $4$-dimensional Riemannian manifold, let $\nabla$ denote its Levi~Civita connection, $R$ its
Riemann tensor and $W$ its Weyl tensor. The latter depends only on the conformal class $[g]$ of the metric $g$
and decomposes as $W=W_++W_-$. Another conformal invariant, the Bach tensor $B$, is the irreducible component
of $\nabla\nabla R$ that, if $M$ is compact, corresponds to the gradient of the Lagrangian \[ g\mapsto
\int\nolimits_M\|W[g]\|^2\,d\upsilon_g,\] in which one may, for topological reasons, replace $W$ by $W_+$ or
$W_-$. The tensor $B$ can be regarded as trace-free symmetric bilinear form, and vanishes if $M$ is self-dual
($W_-=0$) or anti-self-dual ($W_+=0$), in other words if $M$ is half conformally flat. It also vanishes
whenever $[g]$ has an Einstein representative \cite{bach,DER1,besse}. Metrics with zero Bach tensor therefore
form a natural class in which to generalize results on Einstein metrics and curvature flow \cite{TV,CCCMM}.

There are few known examples of `non-trivial' Bach-flat metrics, meaning ones with zero Bach tensor but
satisfying neither of the Weyl and Einstein conditions. A construction of examples is outlined in
\cite{schmidt} and an explicit Lorentzian one is given in \cite{NP}. In the light of the classification of
Einstein, hypercomplex and self-dual metrics on $4$-dimensional Lie groups \cite{Je,BB,barberis,vivian} a
reasonable question is whether there exists a Lie group with such a Bach-flat metric. This note answers the
question by providing two examples. The Lie groups are not unimodular and do not pass to compact
quotients. Nonetheless, our work shows that they have a privileged status in the study of the curvature of
left-invariant metrics in the spirit of \cite{Mi}.

\section{A family of metric Lie algebras}

Our first example generalizes one in \cite{vivian}. Consider the Lie algebra $\ggg_{\alp,\bet}$ defined by a
basis $(e_i)$ of $\R^4$ with non-zero brackets \[[e_2,e_1]=\alp e_2,\quad [e_3,e_1]=\bet e_3,\quad
[e_4,e_1]=(\alp+\bet)e_4,\quad [e_3,e_2]=e_4,\] where $\alp$ and $\bet$ are non-zero real numbers. This gives
a family of solvable, non-nilpotent Lie algebras, which are unimodular if and only if $\alp+\bet=0$. If
$(e^i)$ is the dual basis of $(e_i)$, then one uses the formula $de^i(e_j,e_k)=-e^i([e_j,e_k])$ to encode the
structure constants into the differential system \be\left\{\ba{rcl} de^1 &=& 0\\[-3pt] de^2 &=&\alp e^1\wedge
e^2\\[-3pt] de^3 &=& \bet e^1\wedge e^3\\[-3pt] de^4 &=& (\alp+\bet) e^1\wedge e^4+e^2\wedge e^3.
\ea\right.\ee{esempiosimon} Once one passes to an asociated Lie group, $d$ may be regarded as the exterior
derivative on the space of left-invariant $1$-forms.

\smallbreak

\begin{e-proposition}\label{iso} 
If $\alp\ne \bet$, there exist Lie algebra isomorphisms
\begin{itemize}
\item \ $P_\lambda\colon\ggg_{\alp,\bet}\to\ggg_{\lambda\alp,\lambda \bet}$
  for each $\lambda\ne0$;
\item \ $Q\ \colon\ggg_{\alp,\bet}\to\ggg_{\bet,\alp}$.
\end{itemize}
Any isomorphism $\ggg_{\alp,\bet}\to\ggg_{\alp',\bet'}$ is generated in this way, so $(\alp',\bet')$ equals either
$(\lambda\alp,\lambda \bet)$ or $(\lambda \bet,\lambda\alp)$ for some $\lambda\ne0$.  If $\alp=\bet$ then necessarily
$\lambda=\pm1$.
\end{e-proposition}

\smallbreak

In terms of a basis $(e^i)$ defining $d$, the first isomorphism is given by $P_\lambda(e^1)=e^1/\lambda$ and
$P_\lambda(e^i)=e^i$ for $i>1$. The second, $Q$, merely swaps $e^2$ and $e^3$. To prove the converse, one can
characterize the isomorphism class by the set of $1$-forms satisfying $\omega\wedge d\omega=0$ where
$\omega=\sum_{i=1}^4a_ie^i$.

Endow $\ggg_{\alp,\bet}$ with the inner product for which $(e^i)$ is (dual to) an orthonormal basis.  Two Lie
algebras equipped with inner products can be called \emph{isometric} if there exists a Lie algebra isomorphism
between them that is an isometry. In particular, $\ggg_{\alp,\bet}$ is isometric to both $\ggg_{\bet,\alp}$
and $\ggg_{-\alp,-\bet}$.

\section{Curvature calculations}

Fix non-zero constants $\alp,\bet$, and let $M$ be a Lie group with Lie algebra $\ggg_{\alp,\bet}$. We endow
$M$ with the left-invariant Riemannian metric $g=\sum_{i=1}^4e^{ii}$, where $e^{ii}=e^i\otimes
e^i$. Let $R$ denote its Riemann tensor, and $\rho$ the Ricci tensor, so that
$\rho_{ij}=R^k_{ikj}=R_{kilj}g^{kl}$. From the Cartan structure equations, we obtain
\[\ts\rho=-(\alp^2+\bet^2+\alp\bet)e^{11}-(\alp^2+\alp\bet+\frac14)e^{22}-(\bet^2+\alp\bet+\frac14)e^{33}
-((\alp+\bet)^2-\frac14)e^{44}.\] 
Observe that $g$ is Einstein if and only if $\alp=\beta=\pm\frac12$; this is the symmetric metric on the
complex hyperbolic plane $\C H^2$ \cite{Je}.

The Bach tensor $B$ can be defined by
\[B_{ij}\ =\ \sum_{p=1}^4\nabla_p\nabla_j\rho_{ip} -\frs12\sum_{p=1}^4\nabla_p\nabla_p\rho_{ij}
+\frs13\tau\rho_{ij} -\sum_{p=1}^4\rho_{pi}\rho_{pj} +\frs1{12}\bigg[3\!\sum_{r,s=1}^4
  (\rho_{rs})^2-\tau^2\bigg] \delta_{ij}.\] 
It turns out to be diagonal and, since $\sum_{i=1}^4B_{ii}=0$, we need only record
\[\ba{rcl} 
B_{11} &=&\ts \frac16-\frac16\alp^2+\frac23\alp^3\bet-\frac23\alp^2\bet^2-\frac12\alp\bet+\frac23\alp\bet^3-\frac16\bet^2,\\
B_{22} &=&\ts \frac56\alp^2+\frac12\bet^2+\frac23\alp^3\bet-2\alp\bet^3+\frac76\alp\bet-\frac23\alp^2\bet^2-\frac12,\\
B_{33} &=&\ts \frac12\alp^2+\frac56\bet^2+\frac76\alp\bet-2\alp^3\bet+\frac23\alp\bet^3-\frac23\alp^2\bet^2-\frac12.\\
\ea\]
It is incredibly easy to resolve the system $B=0$, since it yields $(\alp^2-\bet^2)(1+8\alp\bet)=0$. If the second factor is
zero then $\alp,\bet$ are roots of $8x^4-7x^2+\frac18=0$. Let us denote these roots by
$\pm r_1,\pm r_2$ where
\[\ts r_1=\frac14\sqrt{7-3\sqrt5}=\frac18(3\sqrt2-\sqrt{10}\,),\qquad
r_2=-\frac14\sqrt{7+3\sqrt5}=-\frac18(3\sqrt2+\sqrt{10}\,).\]
We obtain eight solutions which, by Proposition~\ref{iso}, fall into three essentially distinct classes:
\be\ts {\rm(i) \ } (\alp,\bet)=\pm(1,1),\qquad
{\rm(ii) \ } (\alp,\bet)=\pm(\frac12,\frac12),\qquad
{\rm(iii) \ } (\alp,\bet)=\pm(r_1,r_2)\hbox{ \ or \ }\pm(r_2,r_1).
\ee{8cases}
Solutions (i) and (ii) correspond to those of \cite{vivian}. The first is hyperhermitian (so $W_+=0$) and the
second is the Einstein metric on $\C H^2$, so the vanishing of their respective Bach tensor is already known.
Our next goal is to show that (iii) is a non-trivial Bach-flat metric.

Because of \cite{vivian}, we know that (iii) is not half conformally flat and this is verified by direct
compution of $W_\pm$. From \cite[Proposition~1]{listing} (see also \cite{KNT}), a necessary condition for a
$4$-dimensional manifold to be locally conformally Einstein is the existence of a non-zero vector field
$T=\sum_{k=1}^4 x_ke_k$ satisfying $(\DIV_4 W)(X,Y,Z)=W(X,Y,Z,T)$, where
\[(\DIV_4 W)(e_h,e_k,e_p) = -\sum_{i=1}^4 \sum_{q=1}^4\bigg[\omega_h^q(e_i)W_{qkpi}+\omega_k^q(e_i)W_{hqpi}
+\omega_p^q(e_i)W_{hkqi}+\omega_i^q(e_i)W_{hkpq}\bigg].\] In case (iii), we discover a contradiction
by examining $(\DIV_4 W)(e_1,e_2,e_j)$ for $j=1$ and $j=2$. Thus,

\smallbreak

\begin{theorem}\label{thm1}  
  Let $G$ be a simply-connected $4$-dimensional Lie group associated to the solvable Lie algebra
  $\ggg_{r_1,r_2}$ defined by {\rm(\ref{esempiosimon})}, and let $h=\sum_{i=1}^4\!e^i\otimes e^i$. The
  Riemannian metric $h$ is non-trivially Bach flat.
\end{theorem}

\smallbreak

\noindent The Ricci tensor of $h$ is diagonal relative to the basis $(e^i)$ with entries $-\frac32$,
$\frac38(-3+\sqrt5)$, $-\frac38(3+\sqrt5)$, $-\frac34$. The tensor $8W_\pm$ has eigenvalues
$2\pm(3\sqrt2-\sqrt{10})$, $2\mp(3\sqrt2+\sqrt{10})$, $-4\pm2\sqrt{10}$.

\section{A 2-step solvable example}\parskip2pt

In order to undertake a more general study of possible left-invariant Bach-flat metrics, one has to work with
an orthonormal basis of $1$-forms satisfying more complicated differential relations.  The solvable case can
be tackled by a case-by-case analysis following Jensen's work on the classification of left-invariant Einstein
metrics \cite{Je}.  The examples (\ref{8cases}) lie in the class of Lie algebras for which $\dim\ggg'=3$ and
$\dim\ggg''=1$ (where $\ggg'=[\ggg,\ggg]$ and $\ggg''=[\ggg',\ggg']$ are the first two terms of the derived
series). Let $G$ be a $4$-dimensional Lie group with such a Lie algebra, admitting a left-invariant Bach-flat
metric $g$.  One can show that there are three possibilities: (i) $g$ is conformally Einstein, (ii) one of
$W_+,W_-$ is zero, or (iii) $\ggg$ is isomorphic to $\ggg_{r_1,r_2}$ and $g$ is homothetic to $h$.

By rotating in the plane $\big<e^2,e^3\big>$ and applying an overall re-scaling, we may convert
$\mathfrak{g}_{r_1,r_2}$ into a form that simplifies the coefficients in $\Q(\sqrt2,\sqrt5)$. Namely, $de^1=0$
and \[ de^2=e^1\wedge e^3,\qquad de^3=e^1\wedge e^2 + \sqrt5 e^1\wedge e^3,\qquad de^4=\sqrt5 e^1\wedge e^4 +
2\sqrt2 e^2\wedge e^3.\] This therefore describes a Bach-flat metric homothetic to $h$. In view of
Theorem~\ref{thm1}, it is natural to ask whether there are solutions in other number fields. In fact, there
exists a 2-step example (meaning that $\ggg''$ vanishes):

\smallbreak

\begin{theorem}\label{thm2} 
  There exists a solvable Lie algebra $\ggg$ defined over $\Q(\sqrt2,\sqrt3)$ with $\dim\ggg'=2$ and
  $\ggg''=\{0\}$, whose associated Lie groups admits a non-trivial Bach-flat metric.
\end{theorem} 

\smallbreak

\noindent The metric Lie algebra in question is given by $de^1=0=de^2$ and
\[\textstyle de^3=\sqrt2 e^1\wedge e^3 + e^2\wedge e^3 -(1/\sqrt3)e^2\wedge e^4,\qquad 
de^4=\sqrt2 e^1\wedge e^4 - e^2\wedge e^4+(1/\sqrt3)e^2\wedge e^3.\]

\smallbreak

The examples of Theorems~\ref{thm1} and \ref{thm2} are obviously restricted. A study of other isomorphism
classes described in \cite{MS} leads us to predict that there are no continuous families of solvable Lie
algebras admitting non-trivial left-invariant Bach-flat metrics, in contrast to the self-dual case
\cite{vivian}.  We can already assert that there are no non-trivial solutions on unimodular Lie groups, since
the unimodular condition simplifies the equations, in particular in the nilpotent cases and in the reductive
cases $\mathfrak{gl}(2,\R)$, $\mathfrak{u}(2)$.

\enddocument